\documentclass{article}
\usepackage{amssymb}
\usepackage{amsfonts}

\input{tcilatex}
\usepackage{polski}
\begin{document}

\begin{center}
\textbf{Homomorpisms and Functional Equations:}

\textbf{The Goldie Equation\\[0pt]
by \\[0pt]
A. J. Ostaszewski}\\[0pt]

\bigskip
\end{center}

\noindent \textbf{Abstract. }The theory of regular variation, in its
Karamata and Bojani\'{c}-Karamata/de Haan forms, is long established and
makes essential use of the Cauchy functional equation. Both forms are
subsumed within the recent theory of Beurling regular variation, developed
elsewhere. Various generalizations of the Cauchy equation, including the Go%
\l \k{a}b-Schinzel functional equation $(GS)$ and Goldie's equation $(GBE)$
below, are prominent there. Here we unify their treatment by
`algebraicization': extensive use of group structures introduced by Popa and
Javor in the 1960s turn all the various (known) solutions into
homomorphisms, in fact identifying them `en passant', and show that $(GS)$
is present everywhere, even if in a thick disguise.

\bigskip


\bigskip

\section{Introduction}

We are concerned here with an `algebraic conversion' of two functional
equations, so that each solution function describes a homomorphism. Both are
known in the functional equations literature in connection originally with
problems arising in utility theory and go back to Lundberg [Lun] and Acz\'{e}%
l [Acz] (see below for more recent studies); there, however, they were
studied in order to classify their solutions, cf. [AczD]. Our purposes here,
which are algebraic, are thus different and are motivated by a different
context. The nearest to our theme of homomorphy is the paper of Kahlig and
Schwaiger [KahS], which studies a sequence of deformations taking the
equation $(GS)$ below in the limit to the classical Cauchy functional
equation $(CFE)$. For us the equations arise first in the classical Karamata
theory of Regular Variation (briefly, RV -- see [BinGT], henceforth BGT, the
standard text, and [BinO3] for updates) and now in the recently developed
theory of Beurling RV, as in [BinO4,6], which includes the Karamata theory.
In the latter context, the first functional equation (in the unknowns $K$
and $\psi ),$ the generalized `Goldie-Beurling equation', has the form%
\begin{equation}
K(x+y\eta (x))-K(y)=\psi (y)K(x)\text{\qquad }(x,y\in \mathbb{R}),
\tag{$GBE$}
\end{equation}%
where for some $\rho \in \mathbb{R}$%
\[
\eta (x)\equiv 1+\rho x,
\]%
the classical Karamata case being $\rho =0$ and the general Beurling case $%
\rho >0.$ In the RV literature this equation appears in [BinG], in work
inspired by Bojani\'{c} and Karamata [BojK], and is due principally to
Goldie (`Goldie's equation'). In both these cases the solution $K$ describes
a function derived from the limiting behaviour of some regularly varying
function (see \S 2 below). For this reason one may expect (by analogy with
the various derivatives encountered in functional analysis) that $K$ should
be an analogue of a linear function. Indeed, for $\rho =0$ and specializing
to $\psi \equiv 1,$ the earliest classical case, $K$ is \textit{additive},
hence a search for homomorphism, when $\rho >0,$ dictates our agenda below
as an algebraic complement (or companion piece) to the analysis of [BinO5,6].

We denote by $GS$ the family of functions $\eta $ above, since they satisfy
the `conditional form ' (as $x,y\geq 0$) of the \textit{Go\l \k{a}%
b-Schintzel functional equation}%
\begin{equation}
\eta (x+y\eta (x))=\eta (x)\eta (y)\text{\qquad }(x,y\in \mathbb{R}_{+}).
\tag{$GS$}
\end{equation}%
For their significance to RV\ see the recent [BinO6] and for their
significance elsewhere, especially to the theory of functional equations,
[Brz5].

The second functional equation of concern replaces the $K$ on the right of $%
(GBE)$ by a further unknown function $\kappa ,$ yielding a natural
`Pexiderized' generalization\footnote{%
Acknowledging the connection, the qualifier $P$ in $(GBE$-$P)$ is for
`Pexiderized' Goldie-Beurling equation -- referring to Pexider's equation: $%
f(xy)=g(x)+h(y)$ and its generalizations -- cf. [Brz1, 3], and the recent
[Jab].}
\begin{equation}
K(x+y\eta (x))-K(y)=\psi (y)\kappa (x)\text{\qquad }(x,y\in \mathbb{R}),
\tag{$GBE$-$P$}
\end{equation}%
considered also in [ChuT]. Passage to this more general format is motivated
by a desire to include ($GS),$ as the case $K=\psi =\eta $ and $\kappa =\eta
-1,$ which turns out to be highly thematic (see [Ost2], and Theorem 1$%
^{\prime }$). Note that in this specialization the corresponding derivatives
are identical: $K^{\prime }=\kappa ^{\prime };$ cf. \S 5.

\section{Popa circle groups}

Recall that any ring $R$ equipped with its \textit{circle product} $x\circ
y:=x+y+xy$ (see [Jac2, II.3], [Jac3]) is a monoid ([Coh1, \S 3]), the
\textit{circle monoid}, with neutral element $0.$ (Below this format is
preferable to the alternative: $x+y-xy,$ isomorphic under the negation $%
x\mapsto -x$.) When $x\circ y=y\circ x=0,$ the elements $x,y$ are \textit{%
quasi-inverses} of each other in the ring, and the \textit{quasi-units}
(those having quasi-inverses) form the \textit{circle group }of $R$. See
e.g. [ColE] for recent advances on circle groups and historical background.
There is an intimate connection with the \textit{Jacobson radical} of a ring
(see [Coh2, \S 5.4], or Jacobson [Jac1] in 1945), characterized as the
maximal ideal of quasi-units, which was prompted by Perlis in 1942
introducing the circle operation. The corresponding notion in a Banach
algebra is that of left and right adverses, similarly defined -- see [Loo,
\S 20C, 21C]. For the connection of adverses with the automatic continuity
of characters, see [Dal, Prop. 3.1].

The \textit{operation}
\[
x\circ _{\eta }y:=x+y\eta (x),
\]%
with $\eta :\mathbb{R}\rightarrow \mathbb{R}$ arbitrary, was introduced in
1965 to the study of equation ($GS$) by Popa [Pop] and later Javor [Jav] (in
the broader context of $\eta :\mathbb{E}\rightarrow \mathbb{F}$, with $%
\mathbb{E}$ a vector space over a commutative field $\mathbb{F}$), who
observed that the equation is equivalent to the operation $\circ _{\eta }$
being associative on $\mathbb{R}$, and that $\circ _{\eta }$ confers a group
structure on $\mathbb{G}_{\eta }:=\{g:\eta (g)\neq 0\}$ -- see [Pop, Prop.
2], [Jav, Lemma 1.2]. Below we term this a \textit{Popa circle group}, or
\textit{Popa group} for short (see \S 3 below), as the case $\eta (x)=1+x$
(i.e. when $\rho =1$ above, so that $\eta $ represents a `shift') yields
precisely the circle group of the ring $\mathbb{R}$.

As $\circ _{\eta }$ turns $\eta $ into a homomorphism from $\mathbb{G}_{\eta
}$ to $(\mathbb{R}\backslash \{0\},\times )$:%
\[
\eta (x\circ _{\eta }y)=\eta (x)\eta (y)\qquad (x,y\in \mathbb{G}_{\eta }),
\]%
and given the group-theoretic framework of RV which leads to $(GBE)$ and $%
(GBE$-$P)$ above, it is natural to seek further group structures permitting $%
(GBE),$ as a property of $K,$ to be \textit{algebraicized }so as to express
a homomorphism between Popa groups:%
\begin{equation}
K(x\circ _{\eta }y)=K(y)\circ _{\sigma }K(x)\text{ for some }\sigma \in GS.
\tag{$CBE$}
\end{equation}%
(Of course, necessarily, $\psi (y)\equiv \sigma (K(y)).)$ Given its context
and form, we term this a \textit{Cauchy-Beurling equation} $(CBE)$. Indeed,
the case $\rho =0$ rewritten as a difference equation,%
\[
\Delta _{y}K(x)-K(y)\psi (x)=0,\text{\qquad }\Delta _{y}K(x):=K(x+y)-K(y),
\]%
suggests that $K(y)$ should induce some form of shift. This difference theme
is exploited in \S 5 on flows, and linked with integration.

Theorem 1 in \S 3 below gives necessary and sufficent conditions for $(GBE)$
to be algebraicized (as above), yielding `en passant' the form of such a $K$
directly from classical results concerning ($CFE$). Likewise Theorem 1$%
^{\prime }$ in \S 4 below gives necessary and sufficent conditions for $GBE$-%
$P$ to be algebraicized (as above); this builds on the technique of Theorem
1, and is similar but more involved. This again yields en passant the form
of such a $K$ directly from classical results concerning ($CFE$).

Theorem 2 in \S 5 continues the investigation of $(GBE$-$P)$ when the
auxiliary $\psi $ is assumed differentiable, reducing it directly to the
context of Theorem 1, i.e. the algebraicization of $(GBE)$. (The
differentiability assumption is again motivated by regular variation.)
Interpreting $\circ _{\eta }$ as a group action, or flow, the underlying
homomorphy is now expressed not by $K$ but by the \textit{relative
flow-velocity} $f(x):=\eta (x)/\psi (x):$ under mild regularity assumptions,
if $K$ solves $(GBE$-$P),$ then $f$ (equivalently its inverse) satisfies%
\[
f(x\circ _{\eta }y)=f(x)f(y).
\]%
There is a converse for $\psi :=\eta /f$ -- see Prop. D in \S 5.

Our quest for algebraicization links with results not only of Acz\'{e}l but
also of Chudziak [Chu1], who in 2006 considered the problem of identifying
pairs $(f,g)$ satisfying the functional equation%
\begin{equation}
f(x+yg(x))=f(x)\circ f(y)  \tag{$ChE$}
\end{equation}%
for $f:\mathbb{R}\rightarrow (S,\circ )$ with $(S,\circ )$ a semigroup, and $%
g:\mathbb{R}\rightarrow \mathbb{R}$ continuous.

We note two recent papers: [Chu2], where $\mathbb{R}$ is replaced by a
vector space over the field of real or complex numbers, and [ChuK], where $%
f,g$ are both assumed continuous and $(S,\circ )$ is the group of
multiplicative reals.

We turn to some background. The functions in $GS$ have their origin for RV\
in the asymptotic analysis of \textit{self-equivarying functions} $\varphi ,$
briefly $\varphi \in SE$ [Ost2], which for some function $\eta >0$ satisfy%
\begin{equation}
\varphi (x+t\varphi (x))/\varphi (x)\rightarrow \eta (t)\text{\qquad }%
(x\rightarrow \infty ,\text{ }\forall t),  \tag{$SE$}
\end{equation}%
locally uniformly in $t$. For $\eta \equiv 1,$ these specialize to the
\textit{self-neglecting} functions of Beurling (BGT 2.3.1, [Kor, IV.11]; cf.
[BinO4]). For $\varphi \in SE$ the limit $\eta =\eta ^{\varphi }$ is
necessarily in $GS$ [Ost2]. Only $(CFE)$ visibly identifies its solution $K$
as a \textit{homomorphism} -- of the additive group $(\mathbb{R},+)$ --
whereas homomorphy is a central feature in the recent topological
development of the theory of regular variation [BinO1,2], [Ost1]. The role
of homomorphy is new in this context, and is one of our principal
contributions here.

At its simplest, a functional equation as above arises when taking limits%
\begin{equation}
K_{F}(t):=\lim\nolimits_{x\rightarrow \infty }[F(x+t\varphi (x))-F(x)],
\tag{$BK$}
\end{equation}%
for $\varphi \in SE$: if $\eta $ satisfies $(SE),$ above then, for $s,t$
ranging over the set $\mathbb{A}$ on which the limit function $K_{F}$ (%
\textit{Beurling kernel}) exists as a locally uniform limit,
\[
K_{F}(s+t)=K_{F}(s/\eta (t))+K_{F}(t):\qquad K_{F}(t+s\eta
(t))=K_{F}(s)+K_{F}(t).
\]%
As we shall see, both $\mathbb{A}$ and $K_{F}(\mathbb{A})$ carry group
structures under which $K_{F}$ is a homomorphism. Thus, even in the
classical context, $(GS)$ plays a significant role albeit disguised and
previously unnoticed, despite its finger-print: the terms +1 or -1,
appearing in the formulas for $K_{F}$ (cf. Th. 1(iv) below). See [BinO2] for
a deeper analysis of the connection between asymptotics of the form $(BK)$
in a general topological setting involving group homomorphisms, and [BinO6]
and [Ost3] for the broader context here.

Previously, in [BinO5], the equations $(GBE$-$P)$ above were all analyzed
using Riemann sums and associated Riemann integrals, introduced there as a
means of extending Goldie's initial approach (via geometric series). Below
we offer an approach to all of the above equations that is new to the
regular variation literature, and partly familiar, albeit in a different
setting, in the $GS$-literature of `addition formulae' -- see [Brz3, 6] and
[Mur] (this goes back to Acz\'{e}l and Go\l \k{a}b [AczG]): we intertwine
Popa groups and integration.

Corresponding to a less restrictive asymptotic analysis (BGT Ch. 3), the
functional equations above give way to functional inequalities. For instance,%
\begin{equation}
F(x+y)\leq e^{y}F(x)+F(y),  \tag{$GFI$}
\end{equation}%
the \textit{Goldie functional inequality} (see [BinO5] for background and
references) becomes group-subadditivity:%
\[
G(x+y)\leq G(x)\circ _{k}G(y).
\]

Our analysis lends new clarification, via the language of homomorphisms, to
the `classical relation' in RV, connecting $K$ and the auxiliary function $%
\psi ,$ which says that $K=c(\psi -1)$ and $\psi \equiv e^{\cdot }$ (cf.
[BGT Lemma 3.2.1], [BinO5, Th. 1]); in particular we point below to the
implicit role of $GS.$ Also, we explain and extend the result of [BinO5, Th.
9] that the solution in $K,$ subject to $K(0)=0,$ positivity of $\kappa $
(i.e. to the right of $0)$ and continuity and positivity of $\psi $,
satisfies for some $c\geq 0$%
\[
K(x)=c\cdot \tau _{f}(x),\text{ for }\tau _{f}(x):=\int_{0}^{x}\mathrm{d}%
u/f(u),\text{ with }f:=\eta /\psi .
\]%
For an interpretation of $\tau _{f}$, inspired by Beck [Bec], as the \textit{%
occupation time measure} (of $[0,x])$ of the continuous $f$-flow: $dx/dt=$ $%
f(x),$ see [BinO6] (and [BinO4]).

\section{Algebraicization of Goldie's equation}

We return to Popa's contribution [Pop], recalling again from Javor [Jav]
that $\circ _{\eta }$ is associative iff $\eta $ satisfies the Go\l \k{a}%
b-Schinzel equation $(GS)$ above. Then for $\eta \neq 0$ $(\mathbb{G}_{\eta
},\circ _{\eta })$ is a group ([Pop, Prop. 2], [Jav, Lemma 1.2]), and $(GS)$
asserts that $\eta $ is a homomorphism from $\mathbb{G}_{\eta }$ to ($%
\mathbb{R}^{\ast },\cdot ):=(\mathbb{R}\backslash \{0\},\times ):$%
\[
\eta (x\circ _{\eta }y)=\eta (x)\eta (y).
\]%
If $\eta $ is injective on $\mathbb{G}_{\eta },$ then $\circ _{\eta }$ is
\textit{commutative}, as $(GS)$ is symmetric on the right-hand side.
Continuous solutions of $(GS),$ positive on $\mathbb{R}_{+},$ are given by
\[
\eta _{\rho }(x):=1+\rho x
\]%
(see e.g. [Brz5] or the more recent [BinO5]). Whenever context permits, if $%
\eta \equiv \eta _{\rho },$ write the group operation and the Popa group as
\[
a\circ _{\rho }b,\qquad (\mathbb{G}_{\rho },\circ _{\rho });
\]%
here $\mathbb{G}_{\rho }=\mathbb{R}\backslash \{1/\rho \}$ and $\mathbb{G}%
_{0}=\mathbb{R}.$ As $(x\circ _{\rho }y)/\rho \rightarrow xy$ as $\rho
\rightarrow \infty ,$ write also $\mathbb{G}_{\infty }:=\mathbb{R}\backslash
\{0\}=\mathbb{R}^{\ast }$, and $\circ _{\infty }\equiv \cdot $
(multiplication)$;$ then $\mathbb{G}_{\rho }$ takes in the additive reals at
one end ($\rho =0),$ and the multiplicative reals at the other; indeed
\[
a\circ _{0}b:=a+b.
\]%
For the intermediate values of $\rho \in (0,\infty ),$ $\eta _{\rho }:%
\mathbb{G}_{\rho }\rightarrow \mathbb{R}^{\ast }$ is an isomorphism, as%
\[
\eta _{\rho }(x\circ _{\rho }y)=\eta _{\rho }(x)\eta _{\rho }(y).
\]%
Rescaling its domain, $\mathbb{G}_{\rho }$ is typified by the case $\rho =1,$
where%
\[
a\circ _{1}b=a+b+ab=(1+a)(1+b)-1:\qquad (\mathbb{G}_{1},\circ _{1})=(\mathbb{%
R}^{\ast },\cdot )-1,
\]%
and the isomorphism $\eta _{1}$ is a shift/translation (cf. [Pop, \S 3]).

Before considering homomorphisms between the groups above we formulate a
result that has two useful variants, relying on \textit{commutativity} or
\textit{associativity}, whence the subscripts.

\bigskip

\noindent \textbf{Lemma 1}$_{\text{com}}$\textbf{.} \textit{If }$(CBE)$\
\textit{holds for an injective }$K,$\textit{\ arbitrary }$\sigma $ \textit{%
with }$\circ _{\sigma }$ \textit{commutative, and arbitrary} $\eta :\mathbb{R%
}\rightarrow \mathbb{R}$, \textit{then} $\eta (u)\equiv 1+\rho u,$ \textit{%
for some constant} $\rho .$

\bigskip

\noindent \textbf{Proof.} Here $K(u+v\eta (u))=K(u)\circ _{\sigma
}K(v)=K(v)\circ _{\sigma }K(u)=K(v+u\eta (v)),$ as $\circ _{\sigma }$ is
commutative. By injectivity, for all $u,v$%
\[
u+v\eta (u)=v+u\eta (v):\quad u(1-\eta (v))=v(1-\eta (u)),
\]%
so $(\eta (u)-1)/u\equiv \rho =$ const. for $u\neq 0,$ and, taking $v=1,$ $%
\eta (u)\equiv 1+\rho u$ for all $u.$ $\square $

\bigskip

\noindent \textbf{Lemma 1}$_{\text{assoc}}$\textbf{.} \textit{If }$(CBE)$\
\textit{holds for an injective }$K,$\textit{\ arbitrary }$\sigma $ \textit{%
with }$\circ _{\sigma }$ \textit{associative, and arbitrary positive
continuous }$\eta :\mathbb{R}\rightarrow \mathbb{R}$, \textit{then} $\eta
(u)=1+\rho u$ $(u\geq 0),$ \textit{for some constant} $\rho .$

\bigskip

\noindent \textbf{Proof. }This follows e.g. from Javor's observation above
connecting associativity with $(GS)$ ([Jav, p. 235]). Here $K(u\circ _{\eta
}(v\circ _{\eta }w))=K(u)\circ _{\sigma }K(v)\circ _{\sigma }K(w)=K((u\circ
_{\eta }v)\circ _{\eta }w),$ so from injectivity:%
\[
u\circ _{\eta }(v\circ _{\eta }w)=(u\circ _{\eta }v)\circ _{\eta }w,
\]%
i.e. $\circ _{\eta }$ is associative, so satisfies $(GS).$ By results in
[Brz2] and [BrzM], positivity and continuity imply $\eta \in GS$. $\square $

\bigskip

For $\circ _{\eta }=\circ _{0}$ and $\circ _{\sigma }=\circ _{\infty },$ the
equation $(CBE)$ reduces to the exponential format of $(CFE)$ ([Kuc, \S %
13.1]; cf. [Jab]). The critical case for Beurling regular variation is for $%
\rho \in (0,\infty ),$ with positive continuous solutions described as
follows. In the table below the four corner formulas correspond to classical
variants of $(CFE).$

\bigskip

\noindent \textbf{Proposition A }(cf. [Chu1])\textbf{. }\textit{For }$\circ
_{\eta }=\circ _{r},\circ _{\sigma }=\circ _{s},$\textit{\ and }$f$ \textit{%
Baire/measurable satisfying }$(CBE),$ \textit{there is }$\gamma \in \mathbb{R%
}$ \textit{so that }$f(t)$ \textit{is given by:}\renewcommand{%
\arraystretch}{1.25}%
\[
\begin{tabular}{|l|l|l|l|}
\hline
Popa parameter & $s=0$ & $s\in (0,\infty )$ & $s=\infty $ \\ \hline
$r=0$ & $\gamma t$ & $(e^{\gamma t}-1)/s$ & $e^{\gamma t}$ \\ \hline
$r\in (0,\infty )$ & $\gamma \log (1+rt)$ & $[(1+rt)^{\gamma }-1]/s$ & $%
(1+rt)^{\gamma }$ \\ \hline
$r=\infty $ & $\gamma \log t$ & $(t^{\gamma }-1)/s$ & $t^{\gamma }$ \\ \hline
\end{tabular}%
\]%
\medskip \newline
\renewcommand{\arraystretch}{1}

\noindent \textbf{Proof. }Each case reduces to $(CFE)$ or a classical
variant by an appropriate shift and rescaling. For instance, given $f,$ for $%
r,s>0$ set%
\[
F(t):=1+sf((t-1)/r):\qquad f(\tau )=(F(1+r\tau )-1)/s.
\]%
Then with $u=1+rx,v=1+ry,$ as $(uv-1)/r=x\circ _{r}y,$%
\[
F(uv)=1+sf(x\circ _{r}y)=1+sf(x)+sf(y)+s^{2}f(x)f(y)=F(u)F(v).
\]%
So, as $F$ is Baire/measurable (see again [Kuc, \S 13]), $F(t)=t^{\gamma }$
and so $f(t)=[(1+rt)^{\gamma }-1]/s.$ The remaining cases are similar. $%
\square $

\bigskip

Th. 1 below is our main result relating solubility of $(GBE)$ and the
existence of a homomorphism, as in condition (iii). Here condition (ii)
identifies the connection $K(u)\equiv (\psi (u)-1)/s$ which is no surprise
in view of [BojK, (2.2)] and BGT Lemma 3.2.1, and the recent [BinO5, Th. 1]
-- cf. \S 1. This, however, is the nub, as $\psi (u)\equiv 1+sK(u)=\eta
_{s}(K(u)).$

Note that (iv) below covers the classical Cauchy case, provided that for $%
\gamma =0$ we interpret both $c(e^{\gamma x}-1)/\gamma $ and $c[(1+\rho
x)^{\gamma }-1]/\rho \gamma $ by continuity as $cx.$

Below and elsewhere a function is \textit{non-trivial} if it not identically
zero and not indentically 1, and is \textit{positive} if it is positive on $%
(0,\infty ).$

\bigskip

\noindent \textbf{Theorem 1. }\textit{For }$\eta \in GS$\textit{\ in the
setting above, }$(GBE)$ \textit{holds for positive }$\psi $\textit{\ and a
non-trivial }$K$ \textit{iff}

\noindent (i) $K$ \textit{is injective;}\newline
\noindent (ii)\textit{\ }$\sigma =:\psi K^{-1}\in GS,$ \textit{equivalently,
either }$\psi \equiv 1,$ \textit{or for some }$s>0$\textit{\ }%
\[
K(u)\equiv (\psi (u)-1)/s\text{ and }\psi (0)=1,\text{ so }K(0)=0;
\]%
\noindent (iii)%
\begin{equation}
K(x\circ _{\eta }y)=K(x)\circ _{\sigma }K(y).  \tag{Hom-1}
\end{equation}%
\textit{Then}\newline
\noindent (iv) \textit{\ for some constants }$c,\gamma $%
\begin{eqnarray*}
\qquad K(x) &\equiv &c\cdot \lbrack (1+\rho x)^{\gamma }-1]/\rho \gamma ,%
\text{ or}\qquad K(t)\equiv \gamma \log (1+\rho t)\qquad (\rho _{\eta }>0),
\\
\text{or}\qquad K(x) &\equiv &c\cdot (e^{\gamma x}-1)/\gamma \qquad (\rho
_{\eta }=0).
\end{eqnarray*}

\bigskip

\noindent \textbf{Proof. }Consider any non-zero $K;$ this is strictly
monotone and so injective, as%
\[
K(x+y)-K(y)=K(x)\psi (y/\eta (x))\qquad (x,y\in \mathbb{G}_{\eta }),
\]%
and so continuous, by [BinO5, Th. 9, or Lemma]. So $\psi $ is continuous,
since%
\[
\psi (y)\equiv \lbrack K(\xi \circ _{\eta }y)-K(y)]/K(\xi ),
\]%
for any $\xi $ with $K(\xi )\neq 0.$ For convenience, write $k:=K^{-1}$ and $%
\sigma (t):=\psi (k(t)),$ i.e. a composition so continuous. Then%
\begin{equation}
K(y)\circ _{\sigma }K(x)=K(y)+\psi (k(K(y))K(x)=K(x)\psi (y)+K(y),  \tag{*}
\end{equation}%
so with $u=K(x),v=K(y),$ $(GBE)$ becomes%
\[
K(k(u)\circ _{\eta }k(v))=v+u\psi (k(v))=v\circ _{\sigma }u:\qquad k(u\circ
_{\sigma }v)=k(u)\circ _{\eta }k(v),
\]%
as $\circ _{\eta }$ is commutative ($\eta \in GS)$. So (Hom-1) follows from
(*). Lemma 1$_{\text{com}}$ now applies to $k$, as $\circ _{\eta }$ is
commutative. So $\sigma \in GS$ (as $\sigma $ is positive and continuous).
So for some $s,\rho \geq 0$%
\[
\eta (t)\equiv 1+\rho t\text{ and }\sigma (t)\equiv 1+st\qquad (t\geq 0).
\]%
That is, $\psi (K^{-1}(t))\equiv \sigma (t)\equiv 1+st,$ so $\psi
(x)=1+sK(x),$ on substituting $t=K(x).$ So if $s>0$%
\[
K(x)=(\psi (x)-1)/s\qquad (x\geq 0).
\]%
If $s=0,$ then $\psi (t)\equiv 1.$ In any case $\psi (0)=1,$ since setting $%
y=0$ in $(GBE)$ gives $(1-\psi (0))K(x)\equiv K(0)=0,$ but $K$ is injective,
so non-trivial. Substituting into ($GBE$) yields (as in [BinO5, Th. 1] for
the case $\circ _{\eta }=+)$
\[
\psi (x\circ _{\eta }y)=\psi (y)(\psi (x)-1)+\psi (y)=\psi (x)\psi (y),
\]%
so $\psi :\mathbb{G}_{\rho }\rightarrow \mathbb{G}_{\infty }$ is a
continuous homomorphism, and Prop. A applies. If $\rho =\rho _{\eta }=0,$
then $\psi (t)\equiv 1$ or $\psi (t)\equiv e^{\gamma t}$ with $\gamma \neq 0,
$ and for $c=\gamma /s$%
\[
K(t)\equiv c(e^{\gamma t}-1)/\gamma ,\qquad (s>0),\text{ or}\qquad
K(t)\equiv \gamma \log (1+\rho t)\qquad (s=0).
\]%
Otherwise, $\psi \equiv (1+\rho x)^{\gamma }$ with $\gamma \neq 0,$ and then
for $c=\rho \gamma /a$%
\[
K(x)\equiv \lbrack (1+\rho x)^{\gamma }-1]/a=c[(1+\rho x)^{\gamma }-1]/\rho
\gamma ,
\]%
with $\gamma =0$ yielding linear $K$ by our `L'Hospital convention'. The
converse is similar but simpler. $\square $

\bigskip

\noindent \textbf{Remarks. }1\textbf{. }For (iv) see [Acz] and [Chu1], and
note from the comparison that all positive solutions arise as homomorphisms.

\noindent 2. Since $0=1_{\mathbb{G}}$ for $\mathbb{G}$ a Popa group, (Hom-1)
implies $K(0)=0.$

The following Goldie functional inequality, for $\eta \in GS$ continuous,
also arises (in Beurling regular variation) for $K:\mathbb{G}_{\eta
}\rightarrow \mathbb{R}$:
\begin{equation}
K(x\circ _{\eta }y)\leq \psi (y)K(x)+K(y)\qquad (x,y\in \mathbb{G}_{\eta })%
\text{;}  \tag{$GBFI$}
\end{equation}%
the case $\eta \equiv 1$ arises in RV (BGT Ch. 3). With $\sigma (x):=\psi
(K^{-1}(x))$ this is%
\[
K(x\circ _{\eta }y)\leq K(x)\circ _{\sigma }K(y)\qquad (x,y\in \mathbb{G}%
_{\eta }).
\]%
Equation $(GFI)$ above has the equivalent form%
\begin{equation}
F(xy)\leq yF(x)+F(y)\qquad (x,y\in \mathbb{R}_{+})\text{,}  \tag{$GFI_{+}$}
\end{equation}%
for $F:\mathbb{R}_{+}\rightarrow \mathbb{R}_{+}.$ The Popa approach with $%
\sigma =F^{-1}$ here yields
\[
F(xy)\leq F(y)\circ _{\sigma }F(x)\qquad (x,y\in \mathbb{R}_{+}),
\]%
i.e. group-theoretic subadditivity (cf. BGT Ch. 3).

\section{Algebraicization of the Pexiderized equation}

In this section, similarly to Th. 1 above (\S 3), we characterize
circumstances when solubility of $(GBE$-$P)$ is equivalent to a homomorphy
(under $K$). This relies on the function $\kappa $ -- see below. A word of
warning: the roles of the functions $\psi $ and $\kappa $ are \textit{%
complementary} rather than equivalent: their interchange forces\ an
interchange of $x$ and $y$ significant for the subtracted term (on the
left). The subsequent section (\S 5) uses $\kappa $ but focuses on $\psi $
(via $\psi /\eta ).$

In $(GBE)$ the value of $\psi (0),$ if non-zero, has no significant role,
and may without loss of generality be scaled to unity; in $(GBE$-$P)$ the
value $\psi (0)$ has a more significant role, both in relation to $\sigma
:=\psi K^{-1}$ (cf. Th. 1) and in controlling whether $K$ and $\kappa $ are
identical. This is clarified by Proposition B below. THrought this section%
\[
\eta (x)\equiv 1+\rho x.
\]%
We begin with some useful

\bigskip

\noindent \textbf{Observations.} 1. \textit{The value of }$K(0)$\textit{\
may be arbitrary; if }$\psi (0)=0,$\textit{\ then }$K$\textit{\ is constant.}

\noindent The solubility of $(GBE$-$P)$ is unaffected by the choice of $%
K(0), $ since $K(x)$ may be replaced by $K(x)-K(0).$ If $\psi (0)=0,$ then $%
K(x)\equiv K(0)$ -- take $y=0$ in $(GBE$-$P).$

\noindent 2. \textit{If }$K$\textit{\ satisfies }$(GBE$\textit{-}$P),$%
\textit{\ then}

\noindent (i) $K(x)\equiv \psi (0)\kappa (x)+K(0);$

\noindent (ii) \textit{provided }$\psi (0)\neq 0,$\textit{\ } $\psi
(0)\kappa $\textit{\ satisfies }$(GBE)$ \textit{and }$\kappa (0)=0.$

\noindent Taking $y=0$ gives (i). Substitution, for $\psi (0)\neq 0,$ into $%
(GBE$-$P)$ yields
\[
\kappa (x+y\eta (x))=\kappa (y)+\kappa (x)\psi (y)/\psi (0).
\]%
In particular, $\kappa (0)=0$ (put $x=y=0).$

\noindent 3. \textit{For }$\psi $\textit{\ and }$\kappa $\textit{\ positive,
both }$K$\textit{\ and }$\kappa $\textit{\ are continuous and invertible.}

\noindent This follows from [BinO5, Lemma] as $\kappa $ here is strictly
monotone; hence so is $K(x)$ by 2(ii). From here we have the following
extension of [BinO5, Th. 1]:

\bigskip

\noindent \textbf{Lemma 2}. \textit{For }$\psi $\textit{\ and }$\kappa $
\textit{positive, there is }$s\geq 0$ \textit{such that }%
\[
\tilde{\psi}(x):=\psi (x)/\psi (0)=1+s\kappa (x).
\]%
\textit{So }$\kappa :\mathbb{G}_{\rho }\rightarrow \mathbb{G}_{s}$\textit{\
is a homomorphism, and }%
\[
\text{either }\tilde{\psi}\equiv 1,\mathit{\ }\text{or }\kappa (x)=(\tilde{%
\psi}(x)-1)/s\text{ with }s>0,
\]%
\textit{equivalently for }$s>0$\textit{, }$\tilde{\psi}:\mathbb{G}_{\rho
}\rightarrow \mathbb{G}_{\infty }$\textit{\ is a homomorphism:}
\[
\tilde{\psi}(x+y\eta (x))=\tilde{\psi}(x)\tilde{\psi}(y).
\]%
\bigskip

\noindent \textbf{Proof. }Since $K(x\circ _{\eta }y)=K(y\circ _{\eta }x)$
and $K(x)=\psi (0)\kappa (x)+K(0),$
\[
\psi (0)\kappa (y)+K(0)+\psi (y)\kappa (x)=\psi (0)\kappa (x)+K(0)+\psi
(x)\kappa (y).
\]%
Since $\kappa $ is positive, for $x,y>0$%
\[
\kappa (y)[1-\tilde{\psi}(x)]=\kappa (x)[1-\tilde{\psi}(y)]:\quad \quad
\lbrack \tilde{\psi}(x)-1]/\kappa (x)=[\tilde{\psi}(y)-1]/\kappa (y)=s,
\]%
say. Substituting $\tilde{\psi}(x)\equiv 1+s\kappa (x)$ in ($GBE$-$P),$%
\[
\kappa (x+y\eta (x))=\kappa (y)+\kappa (x)(1+s\kappa (y)).
\]%
For $s>0,$ writing $\kappa $ in terms of $\tilde{\psi}$ and cancelling $s,$
\[
\lbrack \tilde{\psi}(x+y\eta (x))-1]=[\tilde{\psi}(y)-1]+[\tilde{\psi}(x)-1]%
\tilde{\psi}(y).\qquad \square
\]

\bigskip

\noindent 4. \textit{If }$K=\kappa $\textit{, then }$\psi (0)=1$\textit{\
and }$\sigma (t):=\psi (K^{-1}(t))=1+st.$

\noindent Immediate from 2(ii) above and $\psi (x)/\psi (0)=1+s\kappa (x).$
Theorem 1 in \S 3 above motivates the interest in $\psi K^{-1}$. Proposition
B below extends this observation, and helps clarify Theorem 1$^{\prime },$
the main result of this section.

\bigskip

\noindent \textbf{Proposition B.} \textit{If }$\sigma :=\psi K^{-1}\in GS$%
\textit{\ and }$\psi (t)/\psi (0)=1+s\kappa (t),$\textit{\ then }$\sigma
(t)\equiv 1+st$ \textit{and one of the following two conditions holds:}

\noindent (i)\textit{\ }$\psi (t)\equiv \psi (0)$\textit{\ and }$K(x)\equiv
\psi (0)\kappa (x)+K(0)$\textit{;}

\noindent (ii)\textit{\ }$s>0$ \textit{and }$K=\kappa $\textit{\ iff }$\psi
(0)=1.$

\bigskip

\noindent \textit{Proof.} Put $\sigma (t):=1+ct,$ with $c\geq 0.$ Since $%
\psi (K^{-1}(t))=\psi (0)(1+s\kappa (K^{-1}(t))),$%
\[
1+ct=\psi (0)[1+s\kappa (K^{-1}(t))]:\qquad \lbrack 1+s\kappa (x)]\psi
(0)=(1+cK(x)).
\]%
From the latter, $c=0$ iff $s=0$, as $K$ and $\kappa $ are non-constant. If $%
c=s=0,$ then, again as $K$ is non-constant, $\psi (t)\equiv \psi (0),$ and
so Observation 2(ii) applies. Suppose next that $c>0.$ Then%
\[
K(x)=\psi (0)\kappa (x)s/c+(\psi (0)-1)/c.
\]%
So, again since $K(x)=\psi (0)\kappa (x)+K(0)$ and $\kappa (0)=0,$
\[
K(0)=(\psi (0)-1)/c,\text{ and }s=c>0.
\]%
So if $\psi (0)=1,$ then $K(0)=0$ and $K(x)=\kappa (x).$ Conversely, if $%
K=\kappa ,$ then $K(0)=0,$ and so $\psi (0)=1.$ $\square $

\bigskip

Theorem 1 identifies when $\kappa $ is a homomorphism, but this yields no
similar information about $K$ except when $K(0)=0.$ So we proceed as
follows. Suppose that $(GBE$-$P)$ is soluble with $\psi ,\kappa $ positive.
Then, as $K$ and $\kappa $ are strictly monotone (cf. observation 3 above),
put $y=K^{-1}(v),$ $x=\kappa ^{-1}(u);$ then,%
\[
K(\kappa ^{-1}(u)+K^{-1}(v)\eta (\kappa ^{-1}(u)))=v+\psi
(K^{-1}(v))u=v\circ _{\sigma }u,
\]%
where $\sigma (t):=\psi (K^{-1}(t)).$ Apply $K^{-1}:$
\[
\kappa ^{-1}(u)\circ _{\eta }K^{-1}(v)=\kappa ^{-1}(u)+K^{-1}(v)\eta (\kappa
^{-1}(u))=K^{-1}(v\circ _{\sigma }u).
\]%
Writing%
\[
\alpha (t):=\kappa ^{-1}(K(t)),\qquad \beta (t):=\eta (\kappa
^{-1}(K(t))),\qquad u\equiv K(K^{-1}(u)),
\]%
this says%
\[
\alpha (K^{-1}(u))+K^{-1}(v)\beta (K^{-1}(u))=K^{-1}(v\circ _{\sigma }u).
\]%
This suggests an extension of Popa's idea:%
\[
u\circ v:=\alpha (u)+v\beta (u),
\]%
with $\alpha ,\beta $ continuous and $\alpha $ invertible. Supposing this to
yield a group structure (see below) and assuming $\sigma \in GS$ (so that $%
\circ _{\sigma }$ is commutative), we arrive at a homomorphism%
\begin{equation}
K^{-1}(u\circ _{\sigma }v)=K^{-1}(u)\circ K^{-1}(v).  \tag{Hom-2}
\end{equation}

We need to note the example $\alpha (x)=x+b$ with $\beta (x)\equiv 1.$ Here $%
x\circ y=x+y+b,$ so that $x\circ y\circ z=x+y+z+2b,$ and the neutral element
$e$ satisfies%
\[
x+e+b=x\text{ iff }e=-b,\qquad \text{then }x^{-1}=-x-2b.
\]%
We write $+_{b}$ for this operation and call this group the $b$-\textit{%
shifted additive reals}.\footnote{%
The multiplicative analogue $x\circ y:=xy/b$ comes from the format $x\circ
y:=\alpha (x)+x\beta (y).$} Note that $+_{0}=+=\circ _{0}.$

\bigskip

\noindent \textbf{Proposition C.} \textit{The operation }$\circ $\textit{\
is a group operation on a subset of} $\mathbb{R}$ \textit{containing }$0$%
\textit{\ iff the subset is closed under }$\circ $\textit{\ and for some
constants }$b,c$ \textit{with} $bc=0$
\[
\alpha (x)\equiv x+b\text{ and }\beta (x)\equiv 1+c(x+b).
\]%
\textit{That is:}%
\[
\alpha (x)\equiv x\text{ and }\beta (x)\equiv 1+cx,\text{ OR }\alpha
(x)\equiv x+b\text{ and }\beta (x)\equiv 1.
\]%
\textit{So this is either a Popa group }$x\circ _{c}y:=x+y(1+cx),$ \textit{%
or the }$b$\textit{-shifted additive reals} \textit{with the operation }$%
x+_{b}y:=x+y+b.$

\bigskip

\noindent \textbf{Proof.} Suppose that $\circ $ defines a group. In the
application later we assume that $\alpha $ is injective, but here for $0$ an
element of the group, $\alpha (x)=x\circ 0,$ and then $\alpha $ must be
injective. By associativity,%
\[
(x\circ y)\circ z=\alpha (x\circ y)+z\beta (x\circ y)\qquad (x,y,z\in
\mathbb{R}),
\]%
and%
\[
x\circ (y\circ z)=\alpha (x)+(y\circ z)\beta (x)=\alpha (x)+(\alpha
(y)+z\beta (y))\beta (x)\qquad (x,y,z\in \mathbb{R}).
\]%
Comparing the $z$ terms,%
\[
\beta (x\circ y)=\beta (x)\beta (y)\qquad (x,y\in \mathbb{R}),
\]%
and so%
\begin{equation}
\alpha (x\circ y)=\alpha (x)+\alpha (y)\beta (x)\qquad (x,y\in \mathbb{R}).
\tag{**}
\end{equation}%
So, as $\alpha $ is injective,%
\[
\beta \alpha ^{-1}(\alpha (x)+\alpha (y)\beta (x))=\beta (x)\beta (y)\qquad
(x,y\in \mathbb{R}).
\]%
Put $u:=\alpha (x)$ and $v=\alpha (y):$%
\[
\beta \alpha ^{-1}(u+v\beta \alpha ^{-1}(u))=\beta \alpha ^{-1}(u)\beta
\alpha ^{-1}(v)\qquad (u,v\in \mathbb{R}),
\]%
so that $\beta \alpha ^{-1}\in GS,$ assuming continuity. So for some $c\geq 0
$%
\[
\beta \alpha ^{-1}(u)\equiv 1+cu:\qquad \beta (v)\equiv 1+c\alpha (v)\qquad
(u,v\in \mathbb{R}).
\]%
So%
\[
x\circ y=\alpha (x)+y(1+c\alpha (x))\qquad (x,y\in \mathbb{R}).
\]%
So by (**)%
\[
\alpha (\alpha (x)+y(1+c\alpha (x)))=\alpha (x)+\alpha (y)(1+c\alpha
(x))\qquad (x,y\in \mathbb{R}).
\]%
Recalling that $\beta \alpha ^{-1}(u)\equiv 1+cu,$ writing $u=\alpha (x)$
and $v$ for $y,$ this is%
\[
\alpha (u+v(1+cu))=u+\alpha (v)(1+cu)=u(1+c\alpha (v))+\alpha (v)\qquad
(u,v\in \mathbb{R}).
\]%
Now set $v=0$ to obtain, with $a:=(1+c\alpha (0))$ and $b:=\alpha (0),$%
\[
\alpha (u)=au+b\qquad (u\in \mathbb{R}).
\]%
As $\alpha $ is injective $a\neq 0.$ If $e$ is the neutral element, then%
\[
y=e\circ y=\alpha (e)+y\beta (e)\qquad (y\in \mathbb{R}),
\]%
so $\alpha (e)=0$ (taking $y=0)$ and $\beta (e)=1$ (taking $y\neq 0).$ So $%
\alpha (e)=ae+b=0,$ and so $e=-b/a.$ Right-sided neutrality requires that
\[
x=x\circ e=\alpha (x)+e\beta (x)=ax+b+e(1+cax)=ax-bcx+b+e\qquad (x\in
\mathbb{R}).
\]%
So $e=-b=-b/a,$ so $a=1$ and $bc=0.$

One possibility is $b=0=e$, i.e. $\alpha (x)\equiv x$ and $\beta (x)\equiv
1+cx.$ (Indeed, $e=1_{c}=0.)$ The other possibility is $c=0,$ in which case $%
\beta (x)\equiv 1,$ $\alpha (x)\equiv x+b,$ and $e=-b.$ $\square $

\bigskip

Applying this result we deduce the circumstances when $(GBE$-$P)$ may be
transformed to a homomorphism between (usually, Popa) groups. We then read
off the form of the solution function from Prop. A. In the theorem below we
see that $K(x)\equiv (\psi (y)-1)/s$ only in the cases (i) and (iii), but
not in (ii) -- compare Th. 1. Indeed, in (ii) below $K$ is affinely related
to $\kappa ,$ unless $K(0)=0$ (and then iff $b=0$ and $\kappa \equiv K$).
Section 4 pursues the affine relation.

Below recall that a function is \textit{positive} if it is so on ($0,\infty $%
), and note that in all cases $\kappa $ is a homomorphism between Popa
groups.

\bigskip

\noindent \textbf{Theorem 1}$^{\prime }$ \textit{If} $(GBE$\textit{-}$P)$%
\textit{\ is soluble for }$\psi $ \textit{positive, }$\kappa $\textit{\
positive and invertible, }$\eta (x)\equiv 1+\rho x$\textit{\ (with }$\rho
\geq 0),$\textit{\ then} $\circ $ \textit{is a group operation and }$K^{-1}$%
\textit{\ is a homomorphism under }$\circ $\textit{:}%
\[
K^{-1}(u\circ _{\sigma }v)=K^{-1}(u)\circ K^{-1}(v)\qquad (u,v\in \mathbb{R}%
),
\]%
\textit{\ \noindent iff }$\sigma :=\psi K^{-1}\in GS$\textit{\ and one of
the following three conditions holds:}\newline
\noindent (i) $\rho =0,$\textit{\ }$\circ =\circ _{0}$ \textit{and }$\circ
_{\sigma }=\circ _{s}$ \textit{for some }$s>0$\textit{; then for some }$%
\gamma \in \mathbb{R}$%
\[
K(t)\equiv \kappa (t)\equiv (e^{\gamma t}-1)/s\text{ },\qquad \psi (t)\equiv
e^{\gamma t};
\]

\noindent (ii) $\rho =0,$\textit{\ }$\circ _{\sigma }=\circ _{0}$\textit{and
}$\circ =+_{b}$\textit{for some} $b\in \mathbb{R};$ \textit{then}%
\[
K(t)\equiv \kappa (t+b)=\kappa (t)+\kappa (b),\qquad \psi (t)\equiv 1\qquad
(t\in \mathbb{R}),
\]%
\textit{and }$\kappa :\mathbb{G}_{0}\rightarrow \mathbb{G}_{0}$ \textit{is
linear}$;$\newline
\noindent (iii) $\rho >0,$\textit{\ }$\circ =\circ _{\rho }$ \textit{and }$%
\circ _{\sigma }=\circ _{s}$ \textit{for some }$s\geq 0$\textit{; then for
some }$\gamma \in \mathbb{R}$%
\begin{eqnarray*}
K(t) &\equiv &\kappa (t)\equiv \lbrack (1+\rho t)^{\gamma }-1]/s,\quad
(s>0)\quad ,\text{ or }\gamma \log (1+rt)\quad (s=0), \\
\psi (t) &\equiv &(1+\rho t)^{\gamma }\quad (s>0)\quad ,\text{ or }\psi
(t)\equiv 1\text{ }\quad (s=0).
\end{eqnarray*}

\noindent \textbf{Proof.} We suppose first that $\circ $ is a group
operation. As above%
\[
K^{-1}(v\circ _{\sigma }u)=K^{-1}(u)\circ K^{-1}(v);
\]%
using this and associativity of $\circ $, Lemma 1$_{\text{assoc}}$ (with $%
k=K^{-1}$ for $K,$ $\circ $ for $\circ _{\sigma }$ and $\circ _{\sigma }$
for $\circ _{\eta }$) implies that $\sigma \in GS$, as $\sigma $ is positive
and continuous: so for some $s\geq 0$
\[
\sigma (t)=\psi (K^{-1}(t))=1+st:\qquad \psi (t)=1+sK(t)\quad (t\in \mathbb{R%
}),
\]%
as in Prop. B. So $\circ _{\sigma }$ is commutative and (Hom-2) holds and $%
\psi (0)=1.$

By Prop. C, $K^{-1}$ is a homomorphism iff one of the following two cases
arises.

\noindent \textit{Case (i): Popa case }$\circ =\circ _{c}$\textit{.} For
some $c\geq 0$%
\[
\kappa ^{-1}(K(x))=\alpha (x)\equiv x,\text{ and }\beta (y)=\eta (\kappa
^{-1}(K(y)))\equiv 1+cy.
\]%
So%
\[
K(t)=\kappa (t)\text{ and }1+\rho \kappa ^{-1}(K(t))=1+ct\quad (t\in \mathbb{%
R}).
\]%
If $\rho >0,$ on rearranging $\kappa ^{-1}(K(y))\equiv cy/\rho ,$ so
combining and using injectivity:%
\[
K(t)=\kappa (t)=\kappa (ct/\rho ):\qquad c=\rho \quad (t\in \mathbb{R}).
\]%
So $\circ =\circ _{\rho }$ and by (Hom-2),%
\begin{equation}
K^{-1}(u\circ _{s}v)=K^{-1}(u)\circ _{\rho }K^{-1}(v)\quad (u,v\in \mathbb{R}%
).  \tag{Hom-3}
\end{equation}%
So $K:\mathbb{G}_{\rho }\rightarrow \mathbb{G}_{s}$ is a homomorphism. By
Prop. A for some $\gamma ,$%
\[
K(t)\equiv ((1+\rho t)^{\gamma }-1)/s\quad \text{or }\gamma \log (1+\rho
t)\quad (s=0).
\]%
If $\rho =0,$ then $c=0,$ i.e. $\eta \equiv \beta \equiv 1,$ and so again
(Hom-3) holds but with $\rho =0:$
\[
K(t)=(e^{\gamma t}-1)/s\text{ }\quad (s>0),\text{ or }\gamma t\quad (s=0).
\]

\noindent \textit{Case (ii): Shifted case.}\textbf{\ }For some $b$
\[
\kappa ^{-1}(K(x))=\alpha (x)\equiv x+b,\text{ and }\beta (y)=\eta (\kappa
^{-1}(K(y)))=1+\rho \kappa ^{-1}(K(y))\equiv 1.
\]%
So $\rho =0,$ as $\kappa ^{-1}(K(y))\equiv y+b$ is non-zero. Furthermore, as
$K(x)\equiv \kappa (x+b),$ writing $K$ and $\psi $ in terms of $\kappa $ in $%
(GBE$-$P),$%
\[
K(x+y)=\kappa (x+y+b)=\kappa (y+b)+\kappa (x)(1+s\kappa (y+b))\quad (x,y\in
\mathbb{R}).
\]%
Putting $z=y+b,$
\[
\kappa (x+z)=\kappa (z)+\kappa (x)(1+s\kappa (z))=\kappa (x)\circ _{s}\kappa
(z)\quad (x,z\in \mathbb{R}).
\]%
So\footnote{%
Alternatively, apply Prop. A to $F(t):=K^{-1}(t)+b,$ as $F:\mathbb{G}%
_{s}\rightarrow \mathbb{G}_{0},$ since $K^{-1}(u\circ _{\eta
}v)=K^{-1}(u)+K^{-1}(v)+b.$} $\kappa :\mathbb{G}_{0}\rightarrow \mathbb{G}%
_{s}$ is a homomorphism (and $\kappa (0)=0).$ So if $s=0,$ then $\psi \equiv
1$ and $\kappa $ is linear: $K(x)=\kappa (x+b)=\kappa (x)+\kappa (b).$ If $%
s>0,$ then, as $\psi (0)=1$, $K(x)=\kappa (x)+K(0);$ but by Prop. A%
\[
\kappa (x)=(e^{\gamma x}-1)/s:\qquad K(x)=(e^{\gamma (x+b)}-1)/s=e^{\gamma
b}\kappa (x)+\kappa (b).
\]%
So $b=0,$ and $K=\kappa ,$ which is included as $\circ =\circ _{0}=+_{0}.$
The converse is similar and simpler. $\square $

\bigskip

\noindent \textbf{Remarks.} 1. The implications (i)-(iii) are new here, but
for their conclusions see also [Acz] and [Chu1]; as with Th. 1 in \S 3
above, a comparison shows that all positive solutions arise as homomorphisms.

\noindent 2. The transformations used to obtain a homomorphism in fact
simplify $(GBE$-$P)$ to the case where $\kappa (u)=u$ and $\psi
(v)=1+c\alpha (v)$.

\section{\textbf{Flows}}

Using Riemann sums and their limits [BinO5, Th. 9] gives conditions\footnote{%
Specializing to the present context: $\kappa $ positive to the right near $0$
and $\psi $ continuous.} such that if $(GBE$-$P)$ is soluble, then the
solution function $K$ and the auxiliary $\kappa $ are differentiable, and $%
K^{\prime }=c\cdot \psi /\eta $ for some constant $c.$ We give a new proof
which also extends our understanding of $(GBE$-$P)$ by reference to the
underlying flow velocity $f:=\eta /\psi $.

Indeed, this section focuses via $f$ on the auxiliary $\psi $ rather than on
$\kappa $, though $\kappa $ continues to play a part. We assume below that $%
\psi (0)\neq 0,$ in order to pursue the affine relation $K(x)\equiv \psi
(0)\kappa (x)+K(0)$ (cf. Observation 2(ii) of \S 4). To link results below
to earlier ones note that if $K(0)=0,$ then $K\equiv \kappa $ iff $\psi (0)=1
$ (cf. also Prop. B). Recall, however, that in $(GBE$-$P)$ the value $K(0)$
need not be zero. For $\tau _{f}$ see \S 1.

\bigskip

\noindent \textbf{Theorem 2.} \textit{For }$\mathbb{\kappa },\eta \in GS$
\textit{continuous and }$\psi $ \textit{not identically zero and
differentiable}$:$ \textit{if the solution }$K$\textit{\ to }$(GBE$-$P)$
\textit{is continuous, then either }$K$\textit{\ is constant or:\newline
\noindent }(i)\textit{\ }$K$ \textit{is differentiable and }$K^{\prime
}(x)\equiv \kappa ^{\prime }(0)/f(x)$ \textit{for} $f(x):=\eta (x)/\psi (x)$;%
\newline
\noindent (ii) $\kappa ^{\prime }(x)/\kappa ^{\prime }(0)\equiv K^{\prime
}(x)/K^{\prime }(0)$ \textit{and }$\kappa (x)=c\tau _{f}(x)$\textit{\ for
some }$c\in \mathbb{R}$\textit{;} \newline
\noindent (iii) $K_{0}^{\prime }:=K/K^{\prime }(0):(\mathbb{R},\circ
_{\sigma })\rightarrow (\mathbb{R},\cdot )$\textit{\ is a homomorphism: }%
\[
K_{0}^{\prime }(x+y\eta (x))=K_{0}^{\prime }(x)K_{0}^{\prime }(y).
\]

\bigskip

We defer the proof to the end of the section, but note the immediate

\bigskip

\noindent \textbf{Corollary 1.} \textit{In the setting of Theorem 2\newline
\noindent }(i)\textit{\ }$K(x)\equiv \kappa ^{\prime }(0)\tau _{f}(x)+K(0)$%
\textit{, where }$f(x):=\eta (x)/\psi (x)$ \textit{is the relative
flow-velocity};\newline
\noindent (ii) $\kappa (x)\equiv aK(x)+b$ \textit{for some }$a,b\in \mathbb{R%
};$\newline
\noindent (iii) \textit{provided }$\psi (0)=1$\textit{, the flow-velocity} $%
f:(\mathbb{R},\circ _{\eta })\rightarrow (\mathbb{R},\cdot )$\textit{\ is a
homomorphism, equivalently }$f$\textit{\ solves Chudziak's functional
equation }$(ChE).$ \textit{So }$\psi (x)\equiv \eta (x)/f(x),$ \textit{where
}$f$ \textit{satisfies }$(ChE).$

\bigskip

\noindent \textbf{Remarks. }1. As $(GS)\ $corresponds to $K=\psi =\eta ,$ $%
\mathbb{\kappa }(u)=\eta (u)-1,$ here $K(0)=\psi (0)=\eta (0)=1$ and $f=\eta
/\psi =1;$ so $\tau _{f}(x)=x$ and $\kappa (x)=cx,$ so $\eta (x)=K(x)=\kappa
(x)+K(0)=cx+\eta (0)=1+cx.$

\noindent 2. The classical $(GBE)$ case corresponds to $\eta \equiv 1$ (i.e.
$\rho =0)$ and $\kappa =K,$ so $\psi (x)=1/f(x)=e^{\gamma x}$ , as $f$ is a
Popa-homomorphism by Prop. A. So $K(x)=\kappa (x)=c\tau _{f}(x)$ with $\tau
_{f}(x)\equiv (e^{\gamma x}-1)/\gamma .$

\bigskip

We begin with a Proposition which, taken together with Theorem 2 above,
characterizes the solutions to $(GBE$-$P)$ in terms of $f.$ Below it is more
convenient to take $\mathbb{R}_{+}:=[0,\infty )$.

\bigskip

\noindent \textbf{Proposition D. }\textit{If }$f$ \textit{satisfies }$(CBE)$%
\textit{, then subject to }$K(0)=0$\textit{, }$K\equiv \tau _{f}(x)$\textit{%
\ solves }$(GBE$-$P)$\textit{\ for the auxiliaries }$\psi (x):=\eta (x)/f(x)$%
\textit{\ and }$\mathbb{\kappa }\equiv \tau _{f}(x)$\textit{.}

\bigskip

\noindent \textbf{Proof. }Substituting for $K$ in $(GBE$-$P),$ and using $%
u+\eta \sigma (u)=v+u\eta (v),$ as $\eta \in GS,$ we are to prove that
\[
K(v+u\eta (v))-K(v)=\int_{v}^{v+u\eta (v)}\mathrm{d}t/f(t)=\psi (v)\mathbb{%
\kappa }(u)=\psi (v)\int_{0}^{u}\mathrm{d}t/f(t).
\]%
This follows from%
\begin{eqnarray*}
\psi (v)\int_{0}^{u}\mathrm{d}t/f(t) &=&\eta (v)/f(v)\int_{0}^{u}\mathrm{d}%
t/f(t)=\eta (v)\int_{0}^{u}\mathrm{d}t/f(v)f(t) \\
&=&\eta (v)\int_{0}^{u}\mathrm{d}t/f(v+t\sigma (v))\text{ (put }w=v+t\eta
(v)) \\
&=&\int_{v}^{v+u\eta (v)}\mathrm{d}w/f(w).\qquad \square
\end{eqnarray*}%
\bigskip

\noindent \textbf{Corollary 2. }\textit{In the setting of Prop. D the
solution }$K\equiv \tau _{f}$ \textit{of }$(GBE$-$P)$ \textit{takes one of
the forms:}%
\begin{eqnarray*}
\tau _{f}(x) &\equiv &\int_{0}^{x}e^{\gamma t}\mathrm{d}t=(e^{\gamma
x}-1)/\gamma ,\qquad \qquad \qquad (\rho =0,\gamma \neq 0), \\
\tau _{f}(x) &\equiv &\int_{0}^{x}(1+\rho t)^{\gamma }\mathrm{d}t=((1+\rho
x)^{\gamma +1}-1)/\rho (\gamma +1),\qquad (\rho \in (0,\infty ),\gamma \neq
-1), \\
\tau _{f}(x) &\equiv &x,\qquad (\rho \in \lbrack 0,\infty ]).
\end{eqnarray*}

\noindent \textbf{Proof. }Apply Prop. A, writing $\gamma $ for $-\gamma .$
The final formula is the limit of the cases $\rho =0$ and $\rho >0$ as $%
\gamma $ approaches $0$ or $-1,$ respectively$.$ $\square $

\bigskip

Theorem 2 above is a converse to this. We will need the following
`smoothness result'. (For continuity and differentiability of integrals with
respect to a parameter, see [Jar, \S \S\ 3 and 11]). Recall that for a Popa
group $\mathbb{G}=\mathbb{G}_{\eta },$ $1_{\mathbb{G}}=0$ and $_{\circ }^{-1}
$ denotes its inverse.

\bigskip

\noindent \textbf{Proposition E (Convolution Formula).}\textit{\ For
differentiable }$\eta \in GS,$\textit{\ }%
\[
\mathrm{d}x_{\circ }^{-1}=-\eta (s)^{-2}\mathrm{d}s,
\]%
\textit{so}%
\[
a\ast b(x):=\int_{0}^{x}a(x\circ _{\eta }t_{\circ }^{-1})b(t)\mathrm{d}%
t=\eta (x)\int_{0}^{x}a(s)b(x\circ _{\eta }s_{\circ }^{-1})\frac{\mathrm{d}s%
}{\eta (s)^{2}},
\]%
\textit{for} $a,b$ \textit{continuous; in particular, if }$b$ \textit{is
differentiable/}$\mathcal{C}^{\infty }$\textit{, then so is the convolution
function }$a\ast b,$\textit{\ and}%
\[
a\ast b^{\prime }(x)=\eta ^{\prime }(x)a\ast b(x)+b(0)a(x)/\eta
(x)+\int_{0}^{x}a(s)b^{\prime }(x\circ _{\eta }s_{\circ }^{-1})\frac{\mathrm{%
d}s}{\eta (s)^{3}}.
\]

\bigskip

\noindent \textbf{Proof. }Noting $\eta _{\rho }^{\prime }(x)=\rho ,$
differentiation of $\eta (s_{\circ }^{-1})=1/\eta (s)$ gives%
\[
\rho \mathrm{d}(s_{\circ }^{-1})=-\eta (s)^{-2}\rho \mathrm{d}s.
\]%
Put $s=x\circ t_{\circ }^{-1};$ then $t=x\circ s_{\circ }^{-1}=x+s_{\circ
}^{-1}\eta (x).$ Finally,%
\[
b(x\circ _{\eta }s_{\circ }^{-1})=b(x+s_{\circ }^{-1}\eta (x)),
\]%
which is differentiable in $x;$ so $\mathrm{d}b(x\circ _{\eta }s_{\circ
}^{-1})/\mathrm{d}x=b^{\prime }(x\circ _{\eta }s_{\circ }^{-1})/\eta (s),$
since $1+\rho s_{\circ }^{-1}=\eta (s_{\circ }^{-1})=\eta (s)^{-1}.$ $%
\square $

\bigskip

The next two lemmas prepare the ground for a proof of Theorem 2. The setting
here differs slightly from [BinO5, Th. 9] -- we do not assume non-negativity
of $\kappa ,\psi ,$ but instead freely assume that $\psi $ is
differentiable, since in applications $\psi $ is such (in view of Prop. A).
From this, continuity of $K$ will be shown to imply \textit{automatic
differentiability} -- for a textbook treatment of such matters see J\'{a}rai
[Jar]. We could just as easily assume $\psi $ monotone (also implied by
Prop. A), since a monotone, continuous real function is differentiable
almost everywhere [Rud, \S 8.15] (and is absolutely continuous iff it is the
integral of its derivative).

\bigskip

\noindent \textbf{Lemma 4. }\textit{For continuous }$\mathbb{\kappa },$%
\textit{\ a non-trivial (i.e. non-zero) differentiable function }$\psi ,$
\textit{and continuous }$\eta \in GS:$ \textit{if the solution }$K$\textit{\
to }$(GBE$-$P)$ \textit{is continuous, then }$K$\textit{\ satisfies the
difference equation }%
\[
K(x+u)-K(x)=\mathbb{\kappa }(u/\eta (x))\psi (x),
\]%
\textit{so} $K$\textit{\ has the flow representation}%
\[
xK(x)=\int_{0}^{x}K(t)\mathrm{d}t+\int_{0}^{t}\kappa ((x-t)/\eta (t))\psi (t)%
\mathrm{d}t,
\]%
\textit{\ and so is differentiable on }$\mathbb{R}_{+}.$

\bigskip

\noindent \textbf{Proof. }For $w:=u+v\eta (u)=v+u\eta (v),$ $u=(w-v)/\eta
(v),$ so%
\begin{equation}
K(w)=\mathbb{\kappa }((w-v)/\eta (v))\psi (v)+K(v).  \tag{***}
\end{equation}%
Now write $x$ for $v$ and $u$ for $(w-v)$ to obtain%
\[
K(x+u)-K(x)=\mathbb{\kappa }(u/\eta (x))\psi (x).
\]

In (***) integrate w.r.t. $v$ from $0$ to $w;$ then%
\[
wK(w)=\int_{0}^{w}K(v)\mathrm{d}v+\int_{0}^{w}\mathbb{\kappa }((w-v)/\eta
(v))\psi (v)\mathrm{d}v.
\]%
The second term, being a Beurling convolution, is differentiable by Prop. E.
$\square $

\bigskip

\noindent \textbf{Lemma 5 (Flow Homomorphism).} \textit{If }$K$\textit{\ is
a differentiable solution to }$(GBE$-$P)$\textit{, normalized so that }$%
K(0)=0,$\textit{\ then either }$K\equiv 0,$\textit{\ or:}\newline
\noindent (i)\textit{\ }$K^{\prime }(x)\equiv \kappa ^{\prime }(0)\cdot \psi
(x)/\eta (x)=\kappa ^{\prime }(0)/f(x),$ \textit{for }$f(x)$ \textit{the
flow-velocity }(of \S 1);\newline
\noindent (ii) $K^{\prime }(x)/K^{\prime }(0)\equiv \kappa ^{\prime
}(x)/\kappa ^{\prime }(0),$ \textit{so }$\kappa (x)=c\tau _{f}(x)$ \textit{%
for some }$c\in \mathbb{R};$\newline
\noindent (iii) $K^{\prime }/K^{\prime }(0):(\mathbb{R},\circ _{\sigma
})\rightarrow (\mathbb{R},\cdot )$ \textit{is a homeomorphism:}%
\[
K^{\prime }(x+y\eta (x))=K^{\prime }(x)K^{\prime }(y).
\]%
\textit{In particular, if }$\psi (0)=1,$\textit{\ then\qquad }$f(x+y\eta
(x))=f(x)f(y).$

\bigskip

\noindent \textbf{Proof. }Fixing $y$ with $\psi (y)\neq 0,$ it follows from $%
(GBE$-$P)$ that $\kappa (x)$ is differentiable everywhere. Differentiating
with respect to $x$ and using $x\circ _{\eta }y=y\circ _{\eta }x$
\[
K^{\prime }(x+y\eta (x))\eta (y)=\psi (y)\kappa ^{\prime }(x).
\]%
As $\eta (0)=1,$ substituting 0 alternately for one of $x$ and $y$, and then
for both:
\[
K^{\prime }(y)\eta (y)=\psi (y)\kappa ^{\prime }(0),\qquad K^{\prime
}(x)=\psi (0)\kappa ^{\prime }(x),\qquad K^{\prime }(0)=\psi (0)\kappa
^{\prime }(0).
\]%
So if $\psi (0)\kappa ^{\prime }(0)=0,$ then $K\equiv 0.$ Otherwise,
combining,
\[
\kappa ^{\prime }(x)/\kappa ^{\prime }(0)=K^{\prime }(x)/\kappa ^{\prime
}(0)\psi (0)=K^{\prime }(x)/K^{\prime }(0),
\]%
and in particular $\kappa ^{\prime }(x)=c\psi (x)/\eta (x)=c/f,$ with $%
c=\kappa ^{\prime }(0)/\psi (0).$ So $\kappa (x)=c\tau _{f}(x),$ as $\kappa
(0)=0$ (from $(GBE$-$P)$ for $x=y=0),$ giving (i) and (ii). So
\[
\frac{K^{\prime }(x+y\eta (x))}{K^{\prime }(0)}=\frac{1}{K^{\prime }(0)}%
\frac{\psi (y)}{\eta (y)}\kappa ^{\prime }(x)=\frac{1}{K^{\prime }(0)}\cdot
\frac{K^{\prime }(y)}{\kappa ^{\prime }(0)}\cdot \frac{K^{\prime }(x)}{\psi
(0)}=\frac{K^{\prime }(x)K^{\prime }(y)}{K^{\prime }(0)K^{\prime }(0)},
\]%
equivalently, if $\psi (0)=1$, $K^{\prime }(0)=\kappa ^{\prime }(0)$, $%
f(x)\equiv \kappa ^{\prime }(0)/K^{\prime }(x)$ is a homomorphism$.$ $%
\square $

\bigskip

\noindent \textbf{Proof of Th. 2. }Assuming $K$ non-constant, rescaling if
necessary, without loss of generality $K(0)=0$ and $K^{\prime }(0)=1.$ Now
combine Lemmas 4 and 5. $\square $

\bigskip

\noindent \textbf{Acknowledgements. }We thank a number of colleagues for
constructive comments and references, in particular: Nick Bingham, Charles
Goldie, Mark Roberts, Amol Sasane, Tony Whelan.

\bigskip

\noindent \textbf{References.}\newline
\noindent \lbrack Acz] J. Acz\'{e}l, Extension of a generalized Pexider
equation. \textsl{Proc. Amer. Math. Soc.} \textbf{133} (2005), 3227--3233.%
\newline
\noindent \lbrack AczD] J. Acz\'{e}l, J. Dhombres, \textsl{Functional
equations in several variables. With applications to mathematics,
information theory and to the natural and social sciences.} Encyclopedia of
Math. and its App., 31, CUP, 1989\newline
\noindent \lbrack AczG] J. Acz\'{e}l and S. Go\l \k{a}b, Remarks on
one-parameter subsemigroups of the affine group and their homo- and
isomorphisms, \textsl{Aequat. Math.}, \textbf{4} (1970), 1-10.\newline
\noindent \lbrack Bec] A. Beck, \textsl{Continuous flows on the plane},
Grundl. math. Wiss. \textbf{201}, Springer, 1974. \newline
\noindent \lbrack BinGa] N. H. Bingham, B. Gashi, Logarithmic moving
averages, \textsl{J. Math. Anal. Appl.} \textbf{42} (2015), 1790-1802.%
\newline
\noindent \lbrack BinG] N. H. Bingham, C. M. Goldie, Extensions of regular
variation: I. Uniformity and quantifiers, \textsl{Proc. London Math. Soc.}
(3) \textbf{44} (1982), 473-496.\newline
\noindent \lbrack BinGT] N. H. Bingham, C. M. Goldie and J. L. Teugels,
\textsl{Regular variation}, 2nd ed., Cambridge University Press, 1989 (1st
ed. 1987). \newline
\noindent \lbrack BinO1] N. H. Bingham and A. J. Ostaszewski, The index
theorem of topological regular variation and its applications. \textsl{J.
Math. Anal. Appl.} \textbf{358} (2009), 238-248. \newline
\noindent \lbrack BinO2] N. H. Bingham and A. J. Ostaszewski, Topological
regular variation. I: Slow variation; II: The fundamental theorems; III:
Regular variation. \textsl{Topology and its Applications} \textbf{157}
(2010), 1999-2013, 2014-2023, 2024-2037. \newline
\noindent \lbrack BinO3] N. H. Bingham and A. J. Ostaszewski, The Steinhaus
theorem and regular variation : De Bruijn and after, \textsl{Indagationes
Mathematicae} \textbf{24} (2013), 679-692.\newline
\noindent \lbrack BinO4] N. H. Bingham and A. J. Ostaszewski, Beurling slow
and regular variation, \textsl{Trans. London Math. Soc., }\textbf{1} (2014)
29-56\newline
\noindent \lbrack BinO5] N. H. Bingham and A. J. Ostaszewski, Cauchy's
functional equation and extensions: Goldie's equation and inequality, the Go%
\l \k{a}b-Schinzel equation and Beurling's equation, arxiv.org/abs/1405.3947%
\newline
\noindent \lbrack BinO6] N. H. Bingham and A. J. Ostaszewski, Beurling
moving averages and approximate homomorphisms, arxiv.org/abs/1407.4093.%
\newline
\noindent \lbrack BojK] R. Bojani\'{c} and J. Karamata, \textsl{On a class
of functions of regular asymptotic behavior, }Math. Research Center Tech.
Report 436, Madison, Wis. 1963; reprinted in \textsl{Selected papers of
Jovan Karamata} (ed. V. Mari\'{c}, Zevod Ud\v{z}benika, Beograd, 2009),
545-569.\newline
\noindent \lbrack Brz1] J. Brzd\k{e}k, On the solutions of the functional
equation $f(xf(y)l+yf(x)k)=tf(x)f(y)$, \textsl{Publ. Math. Debrecen},
\textbf{39} (1991), 175-183.\newline
\noindent \lbrack Brz2] J. Brzd\k{e}k, Subgroups of the group $\mathbb{Z}_{n}
$ and a generalization of the Go\l \k{a}b-Schinzel functional equation,
\textsl{Aequat. Math.} \textbf{43} (1992), 59-71.\newline
\noindent \lbrack Brz3] J. Brzd\k{e}k, Some remarks on solutions of the
functional equation $f(x+f(x)^{n}y)=tf(x)f(y),$ \textsl{Publ. Math. Debrecen}%
, \textbf{43} 1-2 (1993), 147-160.\newline
\noindent \lbrack Brz4] J. Brzd\k{e}k, A generalization of the addition
formulae, \textsl{Acta Math. Hungar.,} \textbf{101} (2003), 281-291.\newline
\noindent \lbrack Brz5] J. Brzd\k{e}k, The Go\l \k{a}b-Schinzel equation and
its generalizations, \textsl{Aequat. Math.} \textbf{70} (2005), 14-24.%
\newline
\noindent \lbrack Brz6] J. Brzd\k{e}k, Some remarks on solutions of a
generalization of the addition formulae, \textsl{Aequationes Math.}, \textbf{%
71} (2006), 288--293.\newline
\noindent \lbrack BrzM] J. Brzd\k{e}k and A. Mure\'{n}ko, On a conditional Go%
\l \k{a}b-Schinzel equation, \textsl{Arch. Math.} \textbf{84} (2005),
503-511.\newline
\noindent \lbrack Chu1] J. Chudziak, Semigroup-valued solutions of the Go\l
\k{a}b-Schinzel type functional equation, \textsl{Abh. Math. Sem. Univ.
Hamburg,} \textbf{76} (2006), 91-98.\newline
\noindent \lbrack Chu2] J. Chudziak, Semigroup-valued solutions of some
composite equations, \textsl{Aequat. Math.} \textbf{88} (2014), 183-198.%
\newline
\noindent \lbrack ChuK] J. Chudziak, Z. Ko\v{c}an, Continuous solutions of
conditional composite type functional equations, \textsl{Results Math.}
\textbf{66} (2014), 199--211.\newline
\noindent \lbrack ChuT] J. Chudziak, J. Tabor, Generalized Pexider equation
on a restricted domain. \textsl{J. Math. Psych.} \textbf{52} (2008),
389--392.\newline
\noindent \lbrack Coh1] P. M. Cohn, \textsl{Algebra}, Vol. 1, Weily 1982 (2$%
^{\text{nd}}$ ed.; 1$^{\text{st}}$ ed. 1974).\newline
\noindent \lbrack Coh2] P. M. Cohn, \textsl{Algebra}, Vol. 2, Weily 1989 (2$%
^{\text{nd}}$ ed.; 1$^{\text{st}}$ ed. 1977).\newline
\noindent \lbrack ColE] C. Coleman, D. Easdown, Decomposition of rings under
the circle operation, \textsl{Contributions to Algebra and Geometry} \textbf{%
43} (2002), 55-58.\newline
\noindent \lbrack Dal] H. G. Dales, \textsl{Automatic continuity: a survey},
Bull. London Math. Soc. \textbf{10} (1978),129-183.\newline
\noindent \lbrack Jab] E. Jab\l o\'{n}ska, On solutions of some
generalizations of the Go\l \c{a}b-Schinzel equation. \textsl{Functional
equations in mathematical analysis}, (ed. J. Brzd\k{e}k et al.) 509--521,
Springer, 2012.\newline
\noindent \lbrack Jac 1] N. Jacobson, Radical and Semisimplicity for
arbitrary rings, \textsl{Amer. J. Math.}, 67 (1945), 300-320.\newline
\noindent \lbrack Jac 2] N. Jacobson, \textsl{Lectures in abstract algebra},
Vol. I, Van Nostrand, 1951\newline
\noindent \lbrack Jac 3] N. Jacobson, \textsl{Basic Algebra I}, Freeman, New
York, 1985.\newline
\noindent \lbrack Jar] A. J\'{a}rai, \textsl{Regularity properties of
functional equations in several variables}, Advances in Mathematics 8,
Springer,2005.\newline
\noindent \lbrack Jav] P. Javor, On the general solution of the functional
equation $f(x+yf(x))=f(x)f(y).$ \textsl{Aequat. Math.} \textbf{1} (1968),
235-238.\newline
\noindent \lbrack KahS] P. Kahlig and J. Schwaiger, Transforming the
functional equation of Go\l \c{a}b-Schinzel into one of Cauchy. \textsl{Ann.
Math. Sil.} \textbf{8} (1994), 33-38.\newline
\noindent \lbrack Kor] J. Korevaar, \textsl{Tauberian theorems: A century of
development}. Grundl. math. Wiss. \textbf{329}, Springer, 2004.\newline
\noindent \lbrack Kuc] M. Kuczma, \textsl{An introduction to the theory of
functional equations and inequalities. Cauchy's equation and Jensen's
inequality.} 2nd ed., Birkh\"{a}user, 2009 [1st ed. PWN, Warszawa, 1985].%
\newline
\noindent \lbrack Lun] A. Lundberg, On the functional equation $f(\lambda
(x)+g(y))=%
\mu
(x)+h(x+y)$. \textsl{Aequat. Math.} \textbf{16} (1977), no. 1-2, 21--30.%
\newline
\noindent \lbrack Loo] L. H. Loomis, \textsl{An introduction to abstract
harmonic analysis}, Van Nostrand 1953.\newline
\noindent \lbrack Mur] A. Mure\'{n}ko, On the general solution of a
generalization of the Go\l \k{a}b-Schinzel equation, \textsl{Aequationes
Math.}, \textbf{77} (2009), 107-118.\newline
\noindent \lbrack Ost1] A. J. Ostaszewski, Regular variation, topological
dynamics, and the Uniform Boundedness Theorem, \textsl{Topology Proceedings}%
, \textbf{36} (2010), 305-336.\newline
\noindent \lbrack Ost2] A. J. Ostaszewski, Beurling regular variation, Bloom
dichotomy,and the Go\l \k{a}b-Schinzel functional equation, \textsl{%
Aequationes Math.}, to appear. \newline
\noindent \lbrack Ost3] A. J. Ostaszewski, Asymptotic group actions and
their limits, in preparation.\newline
\noindent \lbrack Pop] C. G. Popa, Sur l'\'{e}quation fonctionelle $%
f[x+yf(x)]=f(x)f(y),$ \textsl{Ann. Polon. Math.} \textbf{17} (1965), 193-198.%
\newline
\noindent \lbrack Rud] W. Rudin, Real and complex analysis, 3rd. ed.,
McGraw-Hill, 1987.

\bigskip

\noindent Mathematics Department, London School of Economics, Houghton
Street, London WC2A 2AE; A.J.Ostaszewski@lse.ac.uk\newpage

\end{document}